# The Frobenius and Monodromy operators for Curves and Abelian Varieties

Robert Coleman, U.C. Berkeley and Adrian Iovita, CICMA Montreal

**Intoduction**

In this paper we will give explicit descriptions of Hyodo and Kato's Frobenius and Monodromy operators on the first $p$-adic de Rham cohomology groups of curves and Abelian varieties with semi-stable reduction over local fields of mixed characteristic. This paper was motivated by the first author's paper [C-pSI], where conjectural definitions of these operators for curves with semi-stable reduction were given. Although B.LeStum wrote a paper entitled "La structure de Hyodo-Kato pour les courbes" ([LS]) he did not prove or claim to have proved that the operators he defined there were the same as Hyodo and Kato's.

The paper is naturally divided into two chapters. In Chapter I, written by the first author, we give the definitions of the Frobenius and Monodromy operators on the de Rham cohomology of Abelian varieties and of curves with semi-stable reduction over a local field $K$.

Suppose for example that $A$ is an Abelian variety with split semi-stable reduction over $K$. If we denote by $K_0$ the maximal unramified subfield of $K$, we define a canonical and functorial $K_0$-lattice in $H^1_{dR}(A)$ denoted $V_0$, and two operators $\Phi_A$ and $N_A$ on $V_0$ such that

i) $\Phi_A$ is $\sigma$-linear, where $\sigma$ is the absolute Frobenius on $K_0$

ii) $N_A$ is $K_0$-linear

and $N_A \Phi_A = p \Phi_A N_A$.

These operators are defined in terms of the "$p$-adic uniformization cross" of $A$

$$\begin{array}{ccccc} & & \Gamma & & \\ & & \downarrow & & \\ T & \to & G & \to & B \\ & & \downarrow & & \\ & & A & & \end{array}$$

where $G$ is a semi-abelian variety, $T$ is a split torus, $B$ is an Abelian variety with good reduction and $\Gamma$ is a free abelian group of finite rank and the diagram makes sense in the rigid analytic category. The monodromy $N_A$ is defined as a residue map along the torus followed by a boundary map and the Frobenius $\Phi_A$ is defined using the Frobenius operators on $T$ and $B$ and the $p$-adic integration of differential forms on $A$. On the other hand, if $X$ is a semi-stable curve over $K$ a Frobenius $\Phi_X$ and a monodromy $N_X$ were defined in [C-pSI]. We prove that if $J$ is the Jacobian of $X$, then $\Phi_X = \Phi_J$ and $N_X = N_J$, where we identify $H^1_{dR}(X)$ and $H^1_{dR}(J)$. In order to prove the identities of these operators one needs to work with de Rham cohomology and duality for 1-motives. These are investigated in section 3.

In Chapter II, which is written by the second author, (it contains essentially the main results of his PhD thesis, Boston University, 1996), the filtered Frobenius Monodromy module attached to $H^1_{dR}(A)$, where $A$ is a split semistable Abelian variety as in Chapter I, is compared to the filtered Frobenius Monodromy module $D_{st}(V(A))^*$ provided by Fontaine's theory, where $V(A) = T_p(A) \otimes_{Z_p} \mathbf{Q}_p$ and $*$ means linear dual. The main result of Chapter II is the following:

The $p$-adic integration pairing

$$<,>: T_p(A) \times H^1_{dR}(A) \to B^+_{dR}$$

defined by P.Colmez in [Cz], induces an isomorphism of filtered Frobenius Monodromy modules between the $K_0$-structure of $H^1_{dR}(A)$ as defined in Chapter I and $D_{st}(V(A))^*$. Our main tool is "the universal covering space" of $A(K)$ defined by

$$\widetilde{A(K)} := \varprojlim(A(\overline{K}), [p]) \times_{A(\overline{K})} A(K)$$

where $[p]$ is the multiplication by $p$ isogeny on $A$. It turns out that $\widetilde{A(K)}_Q := \widetilde{A(K)} \otimes_Z \mathbf{Q}$ is naturally a semistable representation of the Galois group of $\overline{K}$ over $K$ and one can define a map $U\colon T(K) \to D_{st}(\widetilde{A(K)}_Q)$ which plays the role of a "coresidue map" along the torus $T$. In the end we are able to prove, using the results in [C-MP], that Fontaine's monodromy operator on $D_{st}(V(A))$ is essentially induced by Grothendieck's monodromy pairing (after appropriate identifications). As a Corollary we can prove the following

**Theorem.** *Let A be an Abelian variety over the local field K. Then $T_p(A)$ is crystalline if and only if A has good reduction.*

Here $K$ is allowed to be any complete discrete field of characteristic 0 and perfect residue field of characteristic $p$. The only if part of this statement is known by work of J.-M.Fontaine [Fo-BT] and the if part was conjectured by J.-M.Fontaine in [Fo-MGF]. The conjecture was proved in [Fo-MGF] if the ramification index of $K$ is less then $p-1$ and in [M] if $A$ is potentially a product of Jacobians and the residue field of $K$ is finite.

**Aknowledgements** The second author hereby expresses his gratitude towards G.Stevens for his generosity in sharing his ideas and insights with him. During the work on this paper, the second author was visiting U.C.Berkeley and wants to thank this institution for its hospitality.

**Note** Any reference made in one of the chapters to some section or result referes to a section or result in that very chapter except when specifically otherwise stated.

## Chapter I. Definitions of the Operators

In this this section we will recall the definitions of Frobenius and Monodromy operators on the de Rham cohomology groups of curves given in [C-pSI], give definitions of such operators for Abelian varieties and prove that for Jacobians they are equivalent.

Let $K$ be a finite extension of $\mathbf{Q}_p$, $K_0$ the maximal unramified subextension of $K$, $R$ the ring of integers in $K$, $k$ the residue field of $R$ and $v$ the valuation of $K$ which is 1 on a uniformizing parameter of $R$. Suppose $f = [k : \mathbf{F}_p]$. Let $\bar{k}$ be an algebraic closure of $k$ and $\sigma$ the Frobenius automorphism of $\bar{k}/k$. We also use $\sigma$ to denote the lifting of this automorphism to an automorphism of $W(\bar{k})/W(k)$. Also fix a branch log of the $p$-adic logarithm defined over $K$. For a rigid space $S$ over $K$ and a natural number $i$, $H^i_{DR}(S)$ will denote the $i$-th de Rham cohomology group of $S$ over $K$.

### 1. Definitions of $N$ and $F$ for curves

(i) *The monodromy operator*

Suppose $X$ is a connected smooth complete curve over $K$ with a regular semi-stable model $\mathcal{X}$ over $R$ such that the irreducible components of its reduction $\bar{\mathcal{X}}$ are smooth and we will suppose for simplicity of exposition that there are at least two of them and that they, as well as, the singular points of $\bar{\mathcal{X}}$ are defined over $k$. For a subscheme $Y$ of $\bar{\mathcal{X}}$ let $X_Y$ denote the tube of $Y$ considered as a rigid subspace of $X$.

We adopt the notation of Le Stum [LS]. Let $Gr(\mathcal{X})$ be the graph with oriented edges defined as follows: The vertices $V(\mathcal{X})$ of $Gr(\mathcal{X})$ will be the irreducible components of $\bar{\mathcal{X}}$. Let $\bar{\mathcal{X}}^n$ denote the normalization of $\bar{\mathcal{X}}$. Let $m\colon \bar{\mathcal{X}}^n \to \bar{\mathcal{X}}$ be the natural map. The edges $E(X)$ of $Gr(\mathcal{X})$ will be symbols $[x, y]$ where $x$ and $y$ are points on $\bar{\mathcal{X}}^n(\bar{k})$ whose images $\bar{\mathcal{X}}(\bar{k})$ are the same. We set $A([x, y])$ equal to the image of the component of $\bar{\mathcal{X}}^n$ on which $x$ lies and $B([x,y])$ the image in $\bar{\mathcal{X}}$ of the component on which $y$ lies. Then if $e \in E(\mathcal{X})$, $e$ will be an edge ¿from $A(e)$ to $B(e)$. We also define an involution $\tau$ of $E(\mathcal{X})$ by $\tau([x,y]) = [y,x]$.

If $e = [x,y] \in E(\mathcal{X})$ we set $X_e = X_{m(e)}$. We note that $\mathcal{C} = \{X_A \colon A \in V(\mathcal{X})\}$ is an admissible cover of $X_K$ by basic wide opens. We note that since $\mathcal{X}$ is regular, any point in $X(K)$ is contained in a unique element of $\mathcal{C}$.

Let
$$X^0 = \coprod_{A \in V(\mathcal{X})} X_A \quad \text{and} \quad X^1 = \coprod_{e \in E(\mathcal{X})} X_e.$$

Let $\iota$ be the involution on $X^1$, which takes a point in $X_e \subset X^1$ to the corresponding point in $X_{\tau(e)}$. For a module $M$ on which $\iota^*$ acts, $M^{\pm} = \{m \in M \colon \iota^* m = \pm m\}$.

We have a long exact sequence,

$$\to H^0_{DR}(X^0) \xrightarrow{a} H^0_{DR}(X^1)^{-} \xrightarrow{\partial} H^1_{DR}(X) \to H^1_{DR}(X^0) \xrightarrow{b} H^1_{DR}(X^1)^{-} \to \quad (1)$$

For each $e \in E(\mathcal{X})$ we have a natural residue map

$$Res_e \colon H^1_{DR}(X_e) \to H^0_{DR}(X_e).$$

(See [C-RLC].) We set $Res = \bigoplus_{e \in E(\mathcal{X})} Res_e \colon H^1_{DR}(X^1) \to H^0_{DR}(X^1)$. This map takes $H^1_{DR}(X^1)^{+}$ to $H^0_{DR}(X^1)^{-}$.

We define an operatorFrom coleman@math.berkeley.edu Thu Dec 19 19:15:39 1996 Date: Thu, 19 Dec 1996 16:14:50 -0800 (PST) From: "Robert F. Coleman" ¡coleman@math.berkeley.edu¿ To: iovita@math.bu.edu, (Adrian, Iovita) Subject: the file

$N_{\mathcal{X}}$ on $H^1_{DR}(X)$ to be the composition

$$H^1_{DR}(X) \xrightarrow{\rho} H^1_{DR}(X^1)^{+} \xrightarrow{Res} H^0_{DR}(X^1)^{-} \xrightarrow{\partial} H^1_{DR}(C),$$

where $\rho$ is the map obtained from restriction. Let $H(\mathcal{X})$ denote $H^0_{DR}(X^1)^-/a(H^0_{DR}(X^0))$. We will ultimately see that $N_{\mathcal{X}}$ and the image of $H(\mathcal{X})$ in $H^1_{DR}(X)$ are independent of the model $\mathcal{X}$.

(ii) *The Frobenius Operator*

Again we will use the exact sequence (1).

Let $X^\dagger$ denote the dagger completion of $X^0$ along the non-singular locus $NS$ of $\bar{\mathcal{X}}$.

(We note that
$$X_{NS} = \bigcup_{A \in V(\mathcal{X})} (X_A - \bigcup_{A \neq B} X_B)$$
is an underlying affinoid (see [C-RLC]) of $X^0$.) Let $Y$ be a smooth complete curve with a model $\mathcal{Y}$ with good reduction obtained ¿from $X^0$ by glueing in open disks to the ends of $X^0$ (the connected components of $X^0 - X_{NS} \cong X^1$). Then we have a commutative diagram where the rows are exact (for the bottom row see [M, Thm. 4.1]) and the vertical arrows are isomorphisms (see [B1] and [BC]):

$$\begin{array}{ccccccccc}
0 & \to & H^1_{DR}(Y) & \to & H^1_{DR}(X^0) & \to & H^1_{DR}(X^1) & \to & H^2_{DR}(Y) \\
& & \downarrow & & \downarrow & & \downarrow & & \downarrow \\
0 & \to & H^{1\dagger}(Y^\dagger) & \to & H^{1\dagger}(X^\dagger) & \to & K^{E(\mathcal{X})} & \to & K^{V(\mathcal{X})}.
\end{array}$$

Dagger cohomology gives linear Frobenius endomorphisms of the terms in the bottom row (see [MW, Thm. 8.5] and [M, Thm. 4.3]) which is multiplcation by $q$ on $K^{E(\mathcal{X})}$ and $K^{V(\mathcal{X})}$ and is the extension by scalars of the $f$-th power of crystalline Frobenius on

$$H^{1\dagger}(Y^\dagger) \cong H^1_{Cris}(\bar{\mathcal{Y}}, K_0) \otimes_{K_0} K.$$

By virtue of the Riemann hypothesis, the sequence

$$0 \to H^{1\dagger}(Y^\dagger) \to H^{1\dagger}(X^\dagger) \to Ker(K^{E(\mathcal{X})} \to K^{V(\mathcal{X})}) \to 0$$

splits. We can now put a $K_0$/Frobenius structure (by this expression we will mean a $K_0$ sublattice and a $\sigma$-linear Frobenius morphism on it) on $H^1_{DR}(X^0)$ by putting compatible ones on all the terms in the bottom row apart ¿from $H^{1\dagger}(X^\dagger)$. Moreover, $H^1_{DR}(Y)$ is naturally isomorphic to $H^1_{Cris}(\bar{\mathcal{Y}}, K_0) \otimes_{K_0} K$ and

$$K^{E(\mathcal{X})} \cong (W(\bar{k})^{E(\mathcal{X})})^{Gal(L/K_0)} \otimes_{W(k)} K.$$

We actually want and get a $K_0$/Frobenius structure on the kernel of the composition,

$$H^1_{DR}(X^0) \to H^1_{DR}(X^1) \to H^1_{DR}(X^1)^-,$$

which we call $H^1_{DR}(X^0)^+$. In particular, we have an exact sequence

$$0 \to H^1_{DR}(Y) \to H^1_{DR}(X^0)^+ \to H^1_{DR}(X^1)^+ \to H^2_{DR}(Y).$$

**Remark.** *In [C-pSI], we defined another $K_0$/Frobenius structure on $H^1_{DR}(X^0)$ using log-structures. It is probably equivalent but we have not proven that.*

Hence, to get a $K_0$ lattice $V_0$ in $H^1_{DR}(X)$ and a Frobenius operator $F$ on $V_0$ all we have to do is split

$$0 \to H(\mathcal{X}) \to H^1_{DR}(X) \to Ker(H^1_{DR}(X^0) \to H^1_{DR}(X^1)^-) \to 0.$$

We will do this using $p$-adic integration (see [C-pAI] and [C-dS]). Let $\mathcal{W}$ denote the full subcategory of the category of rigid spaces whose objects consist of basic wide opens (see [C-RLC]).

Summarizing results of [C-pAI] and [CdS] we have,

**Theorem 1.1.** *There exists a unique functor ¿from $\mathcal{W}$ to the category of homomorphisms between vector spaces, $W \to \int_W$, where*

$$\int_W : \Omega^1_W(W) \to \mathcal{O}^{loc\,an}_W(W)/K,$$

such that,
$$\int_{\mathbf{A}^1} dz = z \bmod K \quad \text{and} \quad \int_{\mathbf{G}_m} \frac{dz}{z} = \log(z) \bmod K,$$
where $z$ is the standard parameter on the affine line over $K$.

Note: It follows that if $\omega \in W$
$$d \int_W \omega = \omega$$
and if $\omega = df$ where $f$ is rigid analytic, then
$$\int_W df = f \bmod K.$$

If $\omega_A \in \Omega^1_X(X_A)$ for $A \in V(\mathcal{X})$ and $f_e \in \mathcal{O}_X(X_e)$ for $e \in E(\mathcal{X})$,
$$(\{\omega_A\}_{A \in V(\mathcal{X})}, \{f_e\}_{e \in E(\mathcal{X})})$$
will denote the one-hypercochain on $X$ of the complex $\Omega^{\cdot}_X$ with respect to the covering $\mathcal{C}$, $X_A \mapsto \omega_A$ and $(X_A, X_B) \mapsto g_{A,B}$ where $g_{A,B} \in \mathcal{O}_X(X_A \cap X_B)$ is the function such that $g_{A,B}|_{X_e} = f_e$ if $e \in E(\mathcal{X})$ is such that $A(e) = A$ and $B(e) = B$. The hypercochain is a hypercocycle if and only if
$$(\omega_{A(e)} - \omega_{B(e)})|_{X_e} = df_e, \quad \forall e \in E(\mathcal{X}).$$

Let $w \in H^1_{DR}(W)$ and $(\{\omega_A\}, \{f_e\})$ be a one-hypercocycle of $\mathcal{O}_X$ with respect to $\mathcal{C}$ which represents it. We may and will suppose $f_{\tau(e)} = -f_e$. For each $A \in V(\mathcal{X})$, let $s(\omega_A)$ be a representative of $\int_{\log} \omega_A$. Then
$$X_e \mapsto f_e - (s(\omega_{A(e)}) - s(\omega_{B(e)}))$$
represents an element of $H^0_{DR}(X^1)^-$ well defined modulo the image of $H^0_{DR}(X^0)$. Let $I_{\log}(\mathcal{X})$ denote the map which sends $w$ to the image of this class in $H^1_{DR}(X)$. This is the desired splitting.

## 2. $N$ and $F$ for Abelian varieties

(i) *The monodromy operator for Abelian varieties*

Now let $A$ be an Abelian scheme over $K$ with semi-stable reduction. Then we have the "uniformization cross,"

$$\begin{array}{ccccc} & & T & & \\ & & \downarrow & & \\ \Gamma & \to & G & \stackrel{\pi}{\to} & A \\ & & \downarrow & & \\ & & B & & \end{array}$$

where $T$ is a torus, $\Gamma$ is a discrete group, $B$ is an Abelian scheme with good reduction and $G$ is an extension of $B$ by $T$. We have an exact sequence,

$$0 \to Hom(\Gamma, K) \to H^1_{DR}(A) \to H^1_{DR}(G) \to 0. \tag{1}$$

The map from $Hom(\Gamma, K)$ to the kernel of $H^1_{DR}(A) \to H^1_{DR}(G)$ is described as follows: Suppose $\mathcal{C}$ is an admissible covering of $A$ and

$$(\{\omega_U : U \in \mathcal{C}\}, \{f_{UV} : U, V \in \mathcal{C}, U \neq V\})$$

is a one-hypercocycle for $\Omega^{\cdot}_A$ which determines an element of $Ker(H^1_{DR}(A) \to H^1_{DR}(G))$. This means there are functions $h_U$ on $\pi^{-1}U$ for $U \in \mathcal{C}$ such that

$$dh_U = \pi^* \omega_U \quad \text{and} \quad h_U - h_V = \pi^* f_{UV}.$$

Now let $\gamma \in \Gamma$ and $g_U = \gamma^* h_U - h_U$ (this makes sense because $\gamma$ preserves $\pi^{-1}U$). But now,

$$dg_U = 0 \quad \text{and} \quad g_U - g_V = 0.$$

It follows that $\{g_U\}$ corresponds to an element $k_\gamma \in K$. The correspondence $\gamma \mapsto k_\gamma$ is the one we want.

If $h$ is a homomorphism from $\mathbf{G}_m$ into $T$ and $\alpha \in H^1_{DR}(T)$ we set $(\alpha, h) = Res(h^*\alpha)$. This determines an isomorphism from $H^1_{DR}(T)$ onto $Hom_{\mathbf{Z}}(Hom(\mathbf{G}_m, T), K)$. Now there is a perfect pairing

$$Hom(\mathbf{G}_m, T) \times Hom(T, \mathbf{G}_m) \to \mathbf{Z}$$

and so an isomorphism of $Hom_{\mathbf{Z}}(Hom(\mathbf{G}_m, T), K)$ with $Hom(T, \mathbf{G}_m) \otimes_{\mathbf{Z}} K$. Call the isomorphism ¿from $H^1_{DR}(T)$ onto $Hom(T, \mathbf{G}_m) \otimes_{\mathbf{Z}} K$ determined by the above, $Res_T$.

Let $G^0$, $T^0$ and $\mathbf{G}_m^0$ denote the formal completions of $G$, $T$ and $\mathbf{G}_m$ along their special fibers. Then we have an isomorphism

$$f: G(K)/G^0(K) \to T(K)/T^0(K).$$

Moreover, if $h \in Hom(T, \mathbf{G}_m)$, $h$ induces a morphism from $T^0$ to $\mathbf{G}_m^0$. Thus if $v$ is a valuation on $K$, $\gamma \in \Gamma \subset G$ and $h \in Hom(T, \mathbf{G}_m)$ we have an element $(\gamma, h) \in \mathbf{Q}$,

$$v(h(f(\gamma \bmod G^0))).$$

This determines a non-degenerate pairing

$$\Gamma \times Hom(T, \mathbf{G}_m) \to \mathbf{Q},$$

and thus an isomorphism from $Hom(T, \mathbf{G}_m) \otimes K$ onto $Hom(\Gamma, K)$ (see [R-vP]). Finally, let $N_A$ denote the composition

$$H^1_{DR}(A) \to H^1_{DR}(T) \xrightarrow{Res_T} Hom(T, \mathbf{G}_m) \otimes K \to Hom(\Gamma, K) \to H^1_{DR}(A).$$

We note that we have described maps,

$$Hom(\Gamma, K) \to H^1_{DR}(A) \to Hom(T, \mathbf{G}_m) \otimes K. \tag{1}$$

*(ii) The Frobenius operator for Abelian varieties*

We only have to split the exact sequence

$$0 \to Hom(\Gamma, K) \to H^1_{DR}(A) \to H^1_{DR}(G) \to 0. \tag{2}$$

and then put $K_0$/Frobenius structures on $Hom(\Gamma, K)$ and on $H^1_{DR}(G)$.

First we describe the splitting. Suppose $A^\#$ is the universal vectorial extension of $A$ and $G^*$ the pullback of $A^\#$ to $G$. We have,

$$\begin{array}{ccccc}
& & V & = & V \\
& & \downarrow & & \downarrow \\
\Gamma & \to & G^* & \to & A^\# \\
\| & & \downarrow & & \downarrow \\
\Gamma & \to & G & \to & A
\end{array}$$

where $V$ is the vectorial group scheme $V(H^1(\mathcal{O}_A))$ over $K$. If $H$ is a group scheme over a field $L$, we let $Inv_L(H)$ denote the $K$ space of invariant differentials on $H$ over $L$. We have,

$$H^1_{DR}(A) \cong Inv_K(A^\#) \cong Inv_K(G^*).$$

(See Theorem 1.2.2 of [C-DA].)

Using the argument of Bourbaki [III §7.6] one obtains,

**Theorem 2.1.** *If $\omega$ is an invariant differential on $G^*$ over $K$, there is a unique primitive $\lambda_\omega$ of $\omega$ on $G^*(\mathbf{C}_p)$ which is a homomorphism such the restriction of $\lambda_\omega$ to $T$ is contained in*

$$\{\log \circ h \colon h \in Hom(T, \mathbf{G}_m)\} \otimes K.$$

We use this to split (2) as follows: Suppose $\alpha \in H^1_{DR}(A)$ corresponds to the invariant differential $\omega$ on $G^*$. Then $\alpha$ goes to the homomorphism

$$h_\alpha \colon \gamma \in \Gamma \mapsto \lambda_\omega(\gamma).$$

If $\alpha$ is the image of an element $h$ of $Hom(\Gamma, K)$, then $\omega = dg$ where $g \in Hom_K(G^*, \mathbf{G}_a)$ such that $g(\gamma) = h(\gamma)$. It follows that $h_\alpha = h$. Thus the map $I_{\log}(A): \alpha \mapsto h_\alpha$ is a splitting.

Now we let $Hom(\Gamma, \mathbf{Z}) \otimes K \cong Hom(\Gamma, K)$ and for $a \in W(k)$ $\gamma \in \Gamma$, we set

$$F\gamma \otimes a = \gamma \otimes a^\sigma.$$

It remains to determine a $K_0$-structure for $H^1_{DR}(G)$. The schemes $T$, $G$ and $B$ have models with good reduction over $R$. If $X$ is one of these schemes, let $X^\dagger$ denote the dagger completion of $X$ along special fiber of its model. In particular, $B^\dagger = B$. Then we have,

$$\begin{array}{ccccccccc} 0 & \to & H^1_{DR}(B) & \to & H^1_{DR}(G) & \to & H^1_{DR}(T) & \to & 0 \\ & & \downarrow & & \downarrow & & \downarrow & & \\ 0 & \to & H^{1\dagger}(B) & \to & H^{1\dagger}(G^\dagger) & \to & H^{1\dagger}(T^\dagger) & \to & 0 \end{array} \quad (3)$$

We know the top sequence is exact and can check the bottom sequence is as well. Now the outer vertical arrows are isomorphisms. The first is well known, and the last one is easy to check since $T$ is essentially a product of $\mathbf{G}_m$'s (or one can use [BC]). Thus $H^1_{DR}(G)$ is isomorphic to $H^{1\dagger}(G^\dagger)$. Now by Washnitzer-Monsky the objects in the bottom row have compatible actions of Frobenius over $K$. That is, we have endomorphisms $\Phi_B$, $\Phi_G$ and $\Phi_T$ of $H^{1\dagger}(B)$, $H^{1\dagger}(G^\dagger)$ and $H^{1\dagger}(T^\dagger)$ such that the obvious diagrams commute. We now identify the objects on the top row of (3) with the objects directly beneath them. As we will see $\Phi_G$ is a power of the Frobenius operator (tensor $K$) we are after. To make this operator, all we have to do is split the exact sequence (3) since the outer members of this sequence have $W(k)$ structures with $\sigma$-linear Frobenius operators. Suppose $q = |k|$, then $\Phi_T - q$ annihilates $H^1_{DR}(T)$. It follows from the Riemann hypothesis for $B$ that the kernel $M$ of $\Phi_G - q$ in $H^1_{DR}(G)$ maps isomorphically onto $H^1_{DR}(T)$. This gives us the desired splitting.

We saw above that $\Gamma$ maps into $G^*$.

**Proposition 2.2.** $Hom_K(G^*, \mathbf{G}_a) \cong Hom(\Gamma, K)$ under the natural map.

*Proof.* Consider the commutative diagram

$$\begin{array}{ccccccccc}
0 & \to & Hom(\Gamma, K) & \to & H^1_{DR}(J) & \to & H^1_{DR}(G) & \to & 0 \\
& & \downarrow & & \downarrow & & \downarrow & & \\
0 & \to & Hom_K(G^*, \mathbf{G}_a) & \to & Inv_K(G^*) & \to & H^1_{DR}(G^*) & \to & 0
\end{array}$$

in which the rows are exact and the vertical arrows are isomorphisms. It follows $Hom_K(G^*, \mathbf{G}_a) \cong Hom(\Gamma, K)$. The assertion that this map is the natural one follows by chasing the diagram. ∎

## 3. Equality of the Monodromy Operators

Now suppose $X$ is a curve over $K$ with semistable model $\mathcal{X}$ as above and $J$ is the Jacobian of $X$. Then $J$ has semi-stable reduction. Since $H^1_{DR}(X)$ is canonically isomorphic to $H^1_{DR}(J)$ we may consider $N_J$ and $N_\mathcal{X}$ as operators on the same group. We will now show that $N_J = N_\mathcal{X}$.

First let $A$ be an Abelian variety over $K$ with semi-stable reduction. Let the following be the uniformization crosses of $A$ and $\hat{A}$:

$$\begin{array}{ccccc}
& T & & & T' \\
& \downarrow & & & \downarrow \\
\Gamma \to & G & \overset{\pi}{\to} A \quad \text{and} \quad \Gamma' \to & G' & \overset{\pi}{\to} \hat{A}. \\
& \downarrow & & & \downarrow \\
& B & & & \hat{B}
\end{array}$$

Then $\Gamma' = Hom(T, \mathbf{G}_m)$, $T' = Hom(\Gamma, \mathbf{G}_m)$ and $\hat{B}$ is the dual of $B$. We have a canonical pairing $\Gamma \times \Gamma' \to \mathbf{Z}$, called the monodromy pairing which we denote by $(\ ,\ )_{Mon}$. Now $\Gamma'$ is canonically isomorphic to $Hom(T, \mathbf{G}_m)$ which injects onto a lattice of $H^1_{DR}(T)$ via the map which takes $h \in Hom(T, \mathbf{G}_m)$ to the class of $h^*(dT/T)$. We also have a map of $\Gamma$ onto a lattice in $H^1_{DR}(T')$. Thus, by extension of scalars we obtain a pairing $H^1_{DR}(T) \times H^1_{DR}(T')$ with values in $K$. Pulling back via the projections,

$$H^1_{DR}(A) \to H^1_{DR}(T) \quad \text{and} \quad H^1_{DR}(\hat{A}) \to H^1_{DR}(T')$$

we obtain a pairing,
$$H^1_{DR}(A) \times H^1_{DR}(\hat{A}) \to K,$$

which we also call $(\ ,\ )_{Mon}$. Let $(\ ,\ )_{Poin}$ be the cup product (Poincaré) pairing on $H^1_{DR}(A) \times H^1_{DR}(\hat{A})$.

**Theorem 3.1.** *Suppose $\alpha \in H^1_{DR}(A)$ and $\beta \in H^1_{DR}(\hat{A})$. Then*

$$(\alpha, N_{\hat{A}}\beta)_{Poin} = (\alpha, \beta)_{Mon}$$

We also saw that $H^1_{DR}(T) \cong Hom(T, \mathbf{G}_m) \otimes K = \Gamma' \otimes K$. Let $f$ be the natural map from $H^1_{DR}(A)$ into $\Gamma' \otimes K$ and $g$ the map from $Hom(\Gamma', K)$ into $H^1_{DR}(\hat{A})$ as described in section 2. The theorem will follow from result of [C-M], which asserts the the pairing of Raynaud is the same as $(\ ,\ )_{Mon}$ and,

**Lemma 3.2.** *Suppose $\omega \in H^1_{DR}(A)$ and $\rho \in Hom(\Gamma', K)$. Then,*

$$(\omega, g(\rho))_{Poin} = \rho(f(\omega)).$$

*Proof.* A good way to think about this, is in terms of the associated 1-motives, $\Gamma \to G$ and $\Gamma' \to G'$. Let $M =: P\colon X \to H$ be a one motive over $K$. Then Raynaud [R] has defined the de Rham cohomology, $H^1_{DR}(M)$, of $M$ over $K$ as follows: By a vectorial extension of $M$, we mean a commutative diagram

$$\begin{array}{ccc} X & \to & W \\ \| & & \downarrow \\ X & \to & H \end{array}$$

where $W$ is a vectorial extension of $H$. If $M^{\#} =: X \to H^{\#}$ is the universal vectorial extension of $M$, then $H^1_{DR}(M)$ is defined to be the dual of $Lie_K H^{\#}$. Let $\hat{M} =: X' \to H'$ be the 1-motive dual to $M$. We will now define a pairing between the dual spaces of $H^1_{DR}(M)$ and $H^1_{DR}(\hat{M})$. First we explicit $M^{\#}$ and

$\hat{M}^{\#}$. Suppose $M$ is

$$\begin{array}{ccc} & & U \\ & & \downarrow \\ X & \to & H \\ & & \downarrow \\ & & C \end{array}$$

where $U$ is a torus and $C$ is an Abelian variety over $K$. Let $Q\colon X \to C$ be the composition of $P$ with the projection to $C$ and $\hat{C}$ be the dual of $C$. Let $C^{\#}$ be the universal vectorial extension of $C$. Then if, for an extension of an Abelian scheme by a torus, $L$, $\omega_L$ denotes the vectorial scheme whose points over a scheme $S$ consists of the invariant differentials of $L$ defined over $S$, we have an exact sequence

$$0 \to \omega_{\hat{C}} \to C^{\#} \to C \to 0. \tag{1}$$

¿From the exact sequence

$$0 \to U' \to H' \to \hat{C} \to 0,$$

we obtain an exact sequence

$$0 \to \omega_{\hat{C}} \to \omega_{H'} \to \omega_{U'} \to 0.$$

By pushout from (1), we get a vectorial extension

$$0 \to \omega_{H'} \to W \to C \to 0.$$

For a point $Q$ on $C$, let $H'_Q$ denote the corresponding extension of $\hat{C}$ by $\mathbf{G}_m$. A point $\tilde{Q}(x)$, in $C^{\#}(K)$ which maps to $Q(x)$ on $C$, corresponds to a normal invariant differential $\eta_{\tilde{Q}(x)} \in \omega_{H'_{Q(x)}}(K)$, by Theorem 0.3.1 of [C-UVB]. Let $f_x\colon H \to H'_{Q(x)}$ be the map which comes by duality from the map of 1-motives

$$\begin{array}{ccc} \mathbf{Z} & \to & H \\ \downarrow & & || \\ X & \to & H \end{array}$$

where the top arrow is determined by $1 \to Q(x)$ and the left arrow by $1 \mapsto x$. The pair, $(-f_x^* \eta_{\tilde{Q}(x)}, \tilde{Q}(x))$, gives rise to a well defined point $R(x)$ in $W$ which maps to $Q(x)$. By pullback, we get a vectorial extension

$$0 \to \omega_{H'} \to H^{\#} \to H \to 0.$$

Moreover, since $R(x)$ and $P(x)$ both map to $Q(x) \in C$, we get a well defined point $P^{\#}(x) \in H^{\#}$ and $X \to H^{\#}$, $x \mapsto P^{\#}(x) \in H^{\#}$ is the universal vectorial extension of $M$.

Summarizing the above, we have the commutative diagram,

$$\begin{array}{ccccccccc}
0 & \to & \omega_{\hat{C}} & \to & C^{\#} & \to & C & \to & 0 \\
& & \downarrow & & \downarrow & & \downarrow & & \\
0 & \to & \omega_{H'} & \to & W & \to & C & \to & 0 \\
& & \uparrow & & \uparrow & & \uparrow & & \\
0 & \to & \omega_{H'} & \to & H^{\#} & \to & H & \to & 0.
\end{array}$$

¿From this we can identify $Lie_K H^{\#}$ as triples $(a, b, c)$ where $a \in \omega_{H'}(K)$, $b \in Lie_K C^{\#}$ and $c \in Lie_K H$ such that the image of $b$ in $Lie_K C$ equals the image of $c$, modulo the equivalence

$$(a, b, c) \sim (a + d, b - d, c)$$

for $d \in \omega_{\hat{C}}(K)$. We have natural maps,

$$Lie_K(U) \to H^1_{DR}(M) \to \omega_{U'}(K)$$

since $Lie_K(U) \cong Ker(Lie_K(H) \to Lie_K(C))$ and $\omega_{U'}(K) \cong \omega_{H'}(K)/\omega_{\hat{C}}(K)$. Now if $(a, b, c)$ represents an element of $Lie_K(H^{\#})$ and $(a', b', c')$ represents an element in $Lie_K H'^{\#}$, where $X' \to H'^{\#}$ is the universal vectorial extension of $\hat{M}$. We set

$$\langle (a, b, c), (a', b', c') \rangle_M = (a, c') + (b, b') + (c, a'), \qquad (2)$$

where the pairings on the right hand side are the natural ones. This is well defined as
$$(d, c') - (d, b') = 0,$$
for $d \in \omega_{\hat{C}}(K)$, since $b'$ and $c'$ have the same image in $Lie_K \hat{C}$. It is easy to see that $\langle \, , \, \rangle_M$ is non-degenerate.

One can show,

**Theorem 3.3.** *Suppose $M := \Gamma \to G$ is the Raynaud uniformization of the Abelian variety $A$, as above, then there is a natural isomorphism ¿from $H^1_{DR}(A)$ to $H^1_{DR}(M)$ such that (i) the following diagram commutes,*

$$\begin{array}{ccccc} Lie_K(T) & \to & H^1_{DR}(A)\check{\,} & \to & \Gamma \otimes K \\ \| & & \downarrow & & \downarrow \\ Lie_K(T) & \to & H^1_{DR}(M)\check{\,} & \to & \omega_{T'}(K) \end{array}$$

*and (ii) after making the appropriate identifications $\langle \, , \, \rangle_M$ coincides with $(\, , \,)_{Poin}$.*

*Proof.* The isomorphism $H^1_{DR}(A) \to H^1_{DR}(M)$ may be deduced ¿from the observation that the universal vectorial extension of $M$, as an extension of $G$, is the pullback of the universal vectorial extension of $A$. Assertion (i) follows by diagram chasing and one can prove (ii) by reducing to the case of Jacobians where it is not hard to check (although, it would be better to deduce this result by giving a definition of the pairing which is clearly functorial in the analytic category). ∎

This theorem combined with (2) establishes Lemma 3.2 and completes the proof of Theorem 3.1.

We can pull $(\, , \,)_{Mon}$ back to $H^1_{DR}(X)$ and all we have to check to see that the monodromy operators coincide is that, for $\rho$ and $\sigma \in H^1_{DR}(X)$,

$$(\rho, N_\mathcal{X} \sigma)_{Poin} = (\rho, \sigma)_{Mon}. \tag{3}$$

For this, we need a formula for $(\ ,\ )_{Mon}$ on $X$. The image of $H^1_{DR}(X)$ in $H^1_{DR}(X^1)^+$ under the natural map is the kernel of the map to $H^0_{DR}(X^0)$, $(\omega_e)_{e \in E(\mathcal{X})} \mapsto (f_A)$ where

$$f_A = \sum_{\substack{e \in E(\mathcal{X}) \\ A(e)=A}} Res_e \omega_e.$$

Moreover, this kernel is canonically isomorphic to $H_1^{Betti}(Gr(\mathcal{X}), K)$ via

$$(\omega_e) \mapsto {\sum_e}' Res_e(\omega_e) e$$

Here $\sum'$ means the sum over unordered edges ($Res_e(\omega_e)e$ as an element of $H_1^{Betti}(Gr(\mathcal{X}),\ K)$ is independent of the orientation of the edge $e$).

Define a pairing on the free Abelian group on $E(\mathcal{X})$ by setting

$$(e, f) = \begin{cases} 1 & \text{if } e = f \\ -1 & \text{if } f = \tau(e) \\ 0 & \text{otherwise} \end{cases}$$

for edges $e$ and $f$. This induces a pairing on $E^-(X)$ and hence by restriction on $H_1^{Betti}(Gr(X), \mathbf{Z})$. The group $H_1^{Betti}(Gr(X), \mathbf{Z})$ is naturally isomorphic to $\Gamma$ and this is the monodromy pairing [G]. Thus, we obtain a formula for $(\ ,\ )_{Mon}$ on $H^1_{DR}(X)$ using the above mapping of this group into $H_1^{Betti}(Gr(\mathcal{X}),\ K)$ and using the extension by scalars of this formula.

Now suppose $(\{\omega_A\}_{A \in V(\mathcal{X})}, \{f_e\}_{e \in E(\mathcal{X})})$ and $(\{\nu_A\}_{A \in V(\mathcal{X})}, \{g_e\}_{e \in E(\mathcal{X})})$ are hypercocycles representing classes $\omega$ and $\nu$ in $H^1_{DR}(X)$. Then the images of these classes in $H_1^{Betti}(Gr(\mathcal{X}), K)$ are represented by

$${\sum_e}' Res_e(\omega_{A(e)}) e \quad \text{and} \quad {\sum_e}' Res_e(\nu_{A(e)}) e,$$

and so,

$$(\omega, \nu)_{Mon} = {\sum_e}' Res_e(\omega_{A(e)}) Res_e(\nu_{A(e)}).$$

Now, we know $N_\mathcal{X}(\nu)$ is represented by $(\{\alpha_A\}, \{\beta_e\})$ where $\alpha_A = 0$ for all $A$ and $\beta_e = Res_e \nu_{A(e)}$. It follows that

$$(\omega, N_\mathcal{X}(\nu))_{Poin} = {\sum_e}' Res_e(\beta_e \omega_{A(e)})$$
$$= {\sum_e}' Res_e(\omega_{A(e)}) Res_e(\nu_{A(e)}),$$

which completes the proof of (3). It also establishes the equality of the monodromy pairing defined by LeStum in [LS] and that defined in this chapter.

## 4. Equality of Frobenius operators

Suppose now our Abelian variety $A$ is the Jacobian $J$ of a curve $X$ as in section 1. Let $\alpha: X \to J$ be an Albanese morphism.

Suppose $\omega$ is an invariant differential on $G^*$. Then $\omega|_V = dh$ for some homomorphism $h: V \to \mathbf{G}_a$. Now suppose, for $B \in V(\mathcal{X})$, $s_B: X_B \to G^*$ are sections of $G^* \to J$ over $X_B$, i.e., morphisms such that the following diagram commutes,

$$\begin{array}{ccc} X_B & \xrightarrow{s_A} & G^* \\ \downarrow & & \downarrow \\ X & \xrightarrow{\alpha} & J \end{array}.$$

These exist, as one can show that $X_B$ is simply connected, by generalizing the argument of Example 2.5 of [U]. Then for each edge $e \in E(\mathcal{X})$ there exists a unique $\gamma_e \in \Gamma$ such that

$$(s_{A(e)} - s_{B(e)}) X_e \subset V - \gamma_e.$$

For $a \in G^*(K)$, let $T_a$ denote translation by $a$. Let

$$\omega_A = s_A^* \omega \quad \text{and} \quad f_e = \left(T_{\gamma_e} \circ (s_{A(e)} - s_{B(e)})\right)^* h.$$

We know

$$Inv(G^*) \cong Inv(A^\#) \cong H^1_{DR}(J) \cong H^1_{DR}(X).$$

The following proposition makes this isomorphism explicit.

**Proposition 4.1.** *The one-hypercocycle $(\{\omega_A\}, \{f_e\})$ with respect to the covering $\{X_A\}_{A \in V(\mathcal{X})}$ of $X$ represents the de Rham cohomology class on $X$ corresponding to $\omega$.*

Since, as we've stated, the connected components of $X^0$ are simply connected, we see there exists a section $X^0 \to G$ of the diagram

$$\begin{array}{ccc} & & G \\ & & \downarrow \\ X^0 & \to & J \end{array}.$$

Moreover, we can and will assume $X_{NS}$ maps into $G^0$. It follows that we have a commutative diagram,

$$\begin{array}{ccccccccc} 0 & \to & Hom(\Gamma, K) & \to & H^1_{DR}(J) & \to & H^1_{DR}(G) & \to & 0 \\ & & \downarrow & & \downarrow & & \downarrow & & \\ 0 & \to & H(\mathcal{X}) & \to & H^1_{DR}(X) & \to & H^1_{DR}(X^0)^+ & \to & 0 \end{array}$$

Now $H(\mathcal{X})$ is canonically isomorphic to $H^1_{Betti}(Gr(\mathcal{X}), K)$ whose dimension is the rank of $\Gamma$. Since the map from $H^1_{DR}(J)$ to $H^1_{DR}(X)$ is an isomorphism, it follows that all the vertical arrows in this diagram are isomorphisms. We have defined splittings of the rows. We want to show that they are the same. We can make the first vertical arrow more explicit. As we've seen $Hom(\Gamma, K)$ is isomorphic to $Hom_K(G^*, \mathbf{G}_a)$. Suppose $\lambda \in Hom_K(G^*, \mathbf{G}_a)$. Let $\omega = d\lambda$ which is an element of $Inv(G^*)$. Let $\sigma$ denote the hypercocycle made from $\omega$ as above using the sections $s_A$. Then,

$$\sigma - \partial(\{s_A^* \lambda\}),$$

is the hypercocycle $(0, \{a_e\})$ where $a_e = \lambda(\gamma_e)$, where here $\partial$ is the boundary map in hypercohomology.. In summary,

**Lemma 4.2.** *The image of $h \in Hom(\Gamma, K)$ in $H^1_{DR}(X)$ is represented by the hypercocycle $(0, \{h(\gamma_e)\}_{e \in E(\mathcal{X})})$.*

Generalizing Theorem 2.9 of [C-pAI], we have,

**Proposition 4.3.** *Suppose $\omega \in Inv(G^*)$ and $\lambda_\omega$ is the primitive of $\omega$ as specified in Theorem II.1. Then*

$$s_A^* \lambda_\omega + C = \int_{\log} s_A^* \omega.$$

Now suppose $\alpha \in H^1_{DR}(J)$ corresponds to $\omega \in Inv(G^*)$, and $\lambda =: \lambda_\omega$ is its log-primitive. Then $I_{\log}(J)(\alpha)$ is the image of the element in $Hom(\Gamma, K)$, $\gamma \mapsto \lambda(\gamma)$. The image of this in $H^1_{DR}(X)$ is represented by the hypercocycle $(0, \{\lambda(\gamma_e)\})$. On the other hand, as we've seen, the image of $\omega$ in $H^1_{DR}(X)$ is represented by the hypercocycle $(\omega_A, f_e)$, where

$$\omega_A = s_A^* \omega \quad \text{and} \quad f_e = \big(T_{\gamma_e} \circ (s_{A(e)} - s_{B(e)})\big)^* h,$$

$h \in Hom_K(V, \mathbf{G}_a)$ and $dh = \omega|_V$. In fact, by the Theorem 2.1, $h = \lambda|_V$. Then $I_{\log}(\mathcal{X})$ of this class is represented by the hypercocycle which takes vertices to 0 and the edge $e$ to

$$\big(T_{\gamma_e} \circ (s_A - s_B)\big)^* h - (s_A^* \lambda - s_B^* \lambda) = \big(T_{\gamma_e} \circ (s_A - s_B) - (s_A - s_B)\big)^* \lambda$$
$$= \lambda(\gamma_e),$$

where $A = A(e)$ and $B = B(e)$. Thus the two splittings correspond. To conclude, we must show the two $K_0$/Frobenius structures do as well. This is clear for the map $Hom(\Gamma, K) \to H(\mathcal{X})$, since by Lemma 4.2, the image of $Hom(\Gamma, \mathbf{Z})$ in $H(\mathcal{X})$ is fixed by Frobenius.

Since $B$ has good reduction and the connected components of $X^1$ are annuli, the composition

$$H^1_{DR}(B) \to H^1_{DR}(G) \to H^1_{DR}(X^0) \to H^1_{DR}(X^1)$$

is zero. This implies that there is a commutative diagram,

$$\begin{array}{ccccccc}
0 \to & H^1_{DR}(B) & \to & H^1_{DR}(G) & \to & H^1_{DR}(T) & \to 0 \\
& \downarrow & & \downarrow & & \downarrow & \\
0 \to & H^1_{DR}(Y) & \to & H^1_{DR}(X^0)^+ & \to & Ker\big(H^1_{DR}(X^1)^+ \to H^2_{DR}(Y)\big) & \to 0
\end{array}$$

We know the rows are exact and that the map from $H^1_{DR}(G)$ to $H^1_{DR}(X^0)^+$ is an isomorphism. Moreover, using the fact that $X_{NS}$ maps into $G^0$, we see that this map respects the dagger structures we put on these groups. It follows that the splittings we observed of these exact sequences correspond. In particular, the maps ¿from $H^1_{DR}(B)$ to $H^1_{DR}(Y)$ and from $H^1_{DR}(T)$ to $Ker(H^1_{DR}(X^1)^+ \to H^2_{DR}(Y))$ are isomorphisms. The map from $H^1_{DR}(B)$ to $H^1_{DR}(Y)$ respects the crystalline structures since the reduction of $B$ is isomorphic to the product of the Jacobians of the components of the reduction of $Y$ and we can identify this map with the one coming from crystalline cohomology. Also, since $Hom(T, \mathbf{G}_m)$ considered as a subgroup of $H^1_{DR}(T)$ maps into the subgroup of $H^1_{DR}(X^1)^+$ consisting of elements with integral residues we see that this map respects Frobenius as well. Putting this together yields the compatibility of the $K_0$/Frobenius structures we want.